\documentclass[a4,psfig,english,11pt,epsf,portrait]{article}
\usepackage{amsmath,amsfonts,latexsym,amscd,amssymb,theorem}
\usepackage{epsfig}
\usepackage{amsmath}
\usepackage{psfrag}
\def\R{\hbox{\bf R}}
\def\Z{\hbox{\bf Z}}

\def\N{\hbox{\bf N}}

\def\g{\gamma}

\def\O{\Omega}

\def\vphi{\varphi}

\def\<{\langle}
\def\>{\rangle}
\newcommand{\ba}{\begin{eqnarray}}
\newcommand{\ea}{\end{eqnarray}}

\newtheorem{thm}{Theorem}[section]

\newtheorem{theorem}[thm]{Theorem}
\newtheorem{definition}[thm]{Definition}

\newtheorem{rem}[thm]{Remark}

\numberwithin{equation}{section}

\renewcommand{\R}{{\mathbb R}}
\renewcommand{\Z}{{\mathbb Z}}
\renewcommand{\N}{{\mathbb N}}

\renewcommand{\g}{\gamma}

\begin{document}

\title{\bf A remark on a generalization of a logarithmic Sobolev inequality to the H\"{o}lder class}
\author{
\normalsize\textsc{ H. Ibrahim \footnote{ Lebanese University,
Faculty of Sciences, Mathematics Department, Hadeth, Beirut,
Lebanon.
\newline \indent E-mail: ibrahim@cermics.enpc.fr
\newline
\indent $\,\,{}^{1}$LaMA-Liban, Lebanese University, P.O. Box 826
Tripoli, Lebanon.} $^{ ,\,1}$}} \vspace{20pt}

\maketitle

%%%%%%%%%%%%%%%%%%%%%%%%%%%%%%%%%%%%%%%%%%%%%%%%%%%%%%%%%%%%%%%
%
%                  ABSTRACT    ABSTRACT
%
%%%%%%%%%%%%%%%%%%%%%%%%%%%%%%%%%%%%%%%%%%%%%%%%%%%%%%%%%%%%%%%

\centerline{\small{\bf{Abstract}}} \noindent{\small{In a recent work
of the author, a parabolic extension of the elliptic Ogawa type
inequality has been established. This inequality is originated from
the Br\'{e}zis-Gallou\"{e}t-Wainger logarithmic type inequalities
revealing Sobolev embeddings in the critical case. In this paper, we
improve the parabolic version of Ogawa inequality by allowing it to
cover not only the class of functions from Sobolev spaces, but the
wider class of H\"{o}lder continuous functions.}}

\hfill\break
 \noindent{\small{\bf{AMS subject classifications:}}} {\small{42B35,
     54C35, 42B25, 39B05.}}\hfill\break
  \noindent{\small{\bf{Key words:}}} {\small{Littlewood-Paley
      decomposition, logarithmic
      Sobolev inequalities, parabolic $BMO$ spaces, parabolic Lizorkin-Triebel
      spaces, parabolic Besov spaces.}}\hfill\break

%%%%%%%%%%%%%%%%%%%%%%%%%%%%%%%%%%%
%                                 %
%                                 %
%     Section 1                   %
%                                 %
%                                 %
%%%%%%%%%%%%%%%%%%%%%%%%%%%%%%%%%%%

\section{Introduction and main results}\label{sec1}
In \cite{Ibrahim09}, a generalization of the Ogawa type inequality
\cite{Og03} to the parabolic framework has been shown. Ogawa
inequality can be considered as a generalized version in the
Lizorkin-Triebel spaces of the remarkable estimate of
Br\'{e}zis-Gallou\"{e}t-Wainger \cite{BG80, BW80} that holds in a
limiting case of the Sobolev embedding theorem. The inequality
showed in \cite[Theorem 1.1]{Ibrahim09} provides an estimate of the
$L^{\infty}$ norm of a function in terms of its parabolic $BMO$
norm, with the aid of the square root of the logarithmic dependency
of a higher order Sobolev norm. More precisely, for any
vector-valued function $f=\nabla g \in W_{2}^{2m,m}(\R^{n+1})$,
$g\in L^{2}(\R^{n+1})$ with $m,n\in \N^{*}$, $2m>\frac{n+2}{2}$,
there exists a constant $C=C(m,n)>0$ such that:
\begin{equation}\label{Ib:eq1}
\|f\|_{L^{\infty}(\R^{n+1})} \leq C \left(1 + \|f\|_{BMO(\R^{n+1})}
\left(\log^{+} (\|f\|_{W^{2m,m}_{2}(\R^{n+1})} +
\|g\|_{L^{\infty}(\R^{n+1})}) \right)^{1/2} \right),
\end{equation}
where $W^{2m,m}_{2}$ is the parabolic Sobolev space (we refer to
\cite{LSU} for the definition and further properties), and $BMO$ is
the parabolic bounded mean oscillation space (defined via parabolic
balls instead of Euclidean ones \cite[Definition~2.1]{Ibrahim09}). The
above inequality reflects a
limiting case of Sobolev embeddings in the parabolic framework (see
\cite{IJM08, IM09} for similar type inequalities, and \cite{BG80,
BW80, Engler, KOT02, KOT03, KT00, Og03} for various elliptic
versions). By considering functions $f\in W^{2m,m}_{2}(\O_{T})$
defined on the bounded domain
$$\O_{T} = (0,1)^{n}\times (0,T),\quad T>0,$$
we have the following estimate (see \cite[Theorem
1.2]{Ibrahim09}):
\begin{equation}\label{Ib:eq2}
\|f\|_{L^{\infty}(\O_{T})} \leq C \left(1 +
\left(\|f\|_{BMO(\O_{T})} + \|f\|_{L^{1}(\O_{T})}\right)
\left(\log^{+} \|f\|_{W^{2m,m}_{2}(\O_{T})} \right)^{1/2} \right).
\end{equation}
The different norms of $f$ appearing in inequalities (\ref{Ib:eq1})
and (\ref{Ib:eq2}) are finite since
\begin{equation}\label{Ib:eq3}
W^{2m,m}_{2} \hookrightarrow C^{\g,\g/2} \hookrightarrow
L^{\infty}\hookrightarrow BMO \quad\mbox{for some}\quad 0<\g<1,
\end{equation}
where $C^{\g,\g/2}$ is the parabolic H\"{o}lder space that will be
defined later. Moreover, it is easy to check that $g$ bounded and
continuous.

The purpose of this paper to show that the condition $f = \nabla g\in
W^{2m,m}_{2}$ (vector-valued case), or $f\in W^{2m,m}_{2}$
(scalar-valued case) can be relaxed. Indeed, inequalities (\ref{Ib:eq1})
and (\ref{Ib:eq2}) can be applied to a wider class of
H\"{o}lder continuous functions $f = \nabla g\in C^{\g,\g/2}$, $0<\g<1$
(vector-valued case), or $f\in C^{\g,\g/2}$ (scalar-valued case).
To be more precise, we now state the main results of this paper. Our
first theorem is the following:
\begin{theorem}\textit{(Logarithmic H\"{o}lder inequality on
$\R^{n+1}$)}.\label{theo1} Let $0<\g<1$. For any $f = \nabla g \in
C^{\g,\g/2}(\R^{n+1})\cap L^{2}(\R^{n+1})$ with $g\in
L^{2}(\R^{n+1})$, there exists a constant $C = C(\g,n)>0$ such that
\begin{equation}\label{Ib:eq4}
\|f\|_{L^{\infty}(\R^{n+1})} \leq C \left(1 + \|f\|_{BMO(\R^{n+1})}
\left(\log^{+} (\|f\|_{C^{\g,\g/2}(\R^{n+1})} +
\|g\|_{L^{\infty}(\R^{n+1})}) \right)^{1/2} \right).
\end{equation}
\end{theorem}
The second theorem deals with functions defined on the bounded
domain $\O_T$.
\begin{theorem}\textit{(Logarithmic H\"{o}lder inequality on a bounded
    domain)}.\label{theo2} Let
$0<\g<1$.
For any $f \in C^{\g,\g/2}(\O_T)$, there exists a constant $C =
C(\g,n,T)>0$ such that
\begin{equation}\label{Ib:eq5}
\|f\|_{L^{\infty}(\O_T)} \leq C \left(1 + \left(\|f\|_{BMO(\O_T)} +
\|f\|_{L^{1}(\O_T)}\right) \left(\log^{+}
(\|f\|_{C^{\g,\g/2}(\O_T)}) \right)^{1/2} \right).
\end{equation}
\end{theorem}
We notice that inequalities (\ref{Ib:eq4}) and (\ref{Ib:eq5})
directly imply (with the aid of the embeddings (\ref{Ib:eq3}))
(\ref{Ib:eq1}) and (\ref{Ib:eq2}).
\begin{rem}\label{rem1}
The same inequality (\ref{Ib:eq4}) still holds for scalar-valued
functions $f = \frac{\partial g}{\partial x_{i}} \in
C^{\g,\g/2}(\R^{n+1})\cap L^{2}(\R^{n+1})$, $i\in {1, \ldots, n+1}$,
with $g\in L^{\infty}(\R^{n+1})$.
\end{rem}
This paper is organized as follows. In Section \ref{sec2}, we give
the definitions of some basic functional spaces used throughout this
paper. Section \ref{sec3} is devoted to the proofs of the main
results.

%%%%%%%%%%%%%%%%%%%%%%%%%%%%%%%%%%
%
%           SECTION 2
%
%%%%%%%%%%%%%%%%%%%%%%%%%%%%%%%%%%
\section{Definitions}\label{sec2}
Let $\mathcal{O}$ be an open subset of $\R^{n+1}$. A generic element
$z\in \R^{n+1}$ has the form $z=(x,t)$ with $x=(x_{1},\ldots,
x_{n})\in \R^{n}$. We begin by defining parabolic H\"{o}lder spaces
$C^{\g,\g/2}$.
\begin{definition}\textit{(Parabolic H\"{o}lder spaces)}.\label{def1}
For $0<\g<1$, we define the parabolic space of H\"{o}lder continuous
functions of order $\g$ in the following way:
$$C^{\g,\g/2}(\mathcal{O}) = \{f\in C(\overline{\mathcal{O}}),\;
\|f\|_{C^{\g,\g/2}(\mathcal{O})}<\infty\},$$ where
\begin{equation}\label{Hol_norm}
\|f\|_{C^{\g,\g/2}(\mathcal{O})} = \|f\|_{L^{\infty}(\mathcal{O})} +
\<f\>^{(\g)}_{x,\mathcal{O}} + \<f\>^{(\g/2)}_{t,\mathcal{O}},
\end{equation}
with
$$\<f\>^{(\g)}_{x,\mathcal{O}} = \sup_{(x,t), (x',t)\in \mathcal{O},\, x\neq x'} \frac{|f(x,t) -
f(x',t)|}{|x-x'|^{\g}}$$ and
$$\<f\>^{(\g/2)}_{t,\mathcal{O}} = \sup_{(x,t), (x,t')\in \mathcal{O},\, t\neq t'} \frac{|f(x,t) -
f(x,t')|}{|t-t'|^{\g/2}}.$$
\end{definition}
For a detailed study of parabolic H\"{o}lder spaces, we refer the
reader to \cite{LSU}. We now briefly recall some basic facts about
Littlewood-Paley decomposition which are crucial in
obtaining our logarithmic inequalities. Given the expansive
$(n+1)\times(n+1)$ matrix $A = \mbox{diag}\{2,\ldots,2,2^{2}\}$
(parabolic anisotropy), the corresponding Littlewood-Paley
decomposition asserts that any tempered distribution $f\in
\mathcal{S}'(\R^{n+1})$ can be decomposed as
\begin{equation}\label{dyad1}
f=\sum_{j\in \Z} \vphi_{j}*f,\quad \mbox{where}\quad \vphi_{j}(z) =
|\mbox{det} A|^{j} \vphi(A^{j} z),
\end{equation}
with the convergence in $\mathcal{S}'/\mathcal{P}$ (modulo
polynomials). Here $\vphi\in \mathcal{S}(\R^{n+1})$ is a test
function such that $\mbox{supp} \,\hat{\vphi}$ is compact and
bounded away from the origin, and $\sum_{j\in \Z} \hat{\vphi}
(A^{j}z) = 1$ for all $z\in \R^{n+1}\setminus \{0\}$, where
$\hat{\vphi}$ is the Fourier transform of $\vphi$. The sequence
$(\vphi_j)_{j\in \Z}$ is mainly used to define homogeneous
Lizorkin-Triebel and Besov spaces (see for instance \cite{Tri,
Tri1}). However, for defining the inhomogeneous parabolic Besov
space $B^{\g}_{\infty,\infty}$ used later in obtaining our results, we use
a slightly different sequence. Indeed, let $\theta\in
C_{0}^{\infty}(\R^{n+1})$ be any cut-off function satisfying:
\begin{equation}\label{old_neta}
\theta(z) = \left\{
\begin{aligned}
& 1 \quad &\mbox{if}& \quad |z|_{p}\leq 1\\
& 0 \quad &\mbox{if}& \quad |z|_{p}\geq 2,
\end{aligned}
\right.
\end{equation}
where $|\cdot|_{p}$ is the parabolic quasi-norm associated to the
matrix $A$ (see \cite{Ibrahim09}). Taking the new function (but
keeping the same notation) $\vphi_{0}$ defined via the relation
\begin{equation}\label{new_seq}
\hat{\vphi}_{0}= \theta,
\end{equation}
we can give the definition of the Besov space
$B^{\g}_{\infty,\infty}$.
\begin{definition}\textit{(Parabolic inhomogeneous Besov
spaces)}.\label{def2} Take the smoothness parameter $0<\g<1$. Let
$(\vphi_{j})_{j\in \Z}$ be the sequence such that $\vphi_{0}$ is
given by (\ref{new_seq}), while $\vphi_{j}$ is given by
(\ref{dyad1}) for all $j\geq 1$. We define the parabolic
inhomogeneous Besov space $B^{\g}_{\infty,\infty}$ as the space of
all functions $f\in \mathcal{S}'(\R^{n+1})$ with finite quasi-norms
$$
\|f\|_{B^{\g}_{\infty,\infty}} = \sup_{j\geq 0} 2^{\g j} \|\vphi_{j}
* f\|_{L^{\infty}(\R^{n+1})}.
$$
\end{definition}

%%%%%%%%%%%%%%%%%%%%%%%%%%%%%%%%%%
%
%           SECTION 3
%
%%%%%%%%%%%%%%%%%%%%%%%%%%%%%%%%%%
\section{Proofs of theorems}\label{sec3}
We begin with the proof of Theorem~\ref{theo1} that strongly relies
on the results obtained in \cite{Ibrahim09}.\\

\noindent \textbf{Proof of Theorem~\ref{theo1}.} Let $N\in \N$ be
any arbitrary integer. Using (\ref{dyad1}), we estimate
$\|f\|_{L^{\infty}}$ in the following way:
\begin{eqnarray}\label{fer0}
\|f\|_{L^{\infty}} &\leq& \Big\|\sum_{j<-N}2^{\g j}2^{-\g
j}|\vphi_{j}
*
  f|\Big\|_{L^{\infty}} +\Big\|\sum_{|j|\leq N}|\vphi_{j} *
  f|\Big\|_{L^{\infty}} + \Big\|\sum_{j>N}2^{-\g j}2^{\g j}|\vphi_{j} *
  f|\Big\|_{L^{\infty}}\nonumber
\\
&\leq& C_{\g}2^{-\g N}\overbrace{\Big\|\Big(\sum_{j<-N} 2^{-2\g
j}|\vphi_{j}
*
  f|^{2}\Big)^{1/2}\Big\|_{L^{\infty}}}^{A_1} + (2N+1)^{1/2}
  \overbrace{\Big\|\Big(\sum_{|j|\leq N}|\vphi_{j} *
  f|^{2}\Big)^{1/2}\Big\|_{L^{\infty}}}^{A_2}\nonumber\\
& &+ \,C'_{\g} 2^{-\g N} \overbrace{\big(\sup_{j>N} 2^{\g j}
\|\vphi_{j}
* f\|_{L^{\infty}}\big)}^{A_3},
\end{eqnarray}
where
$$ C_{\g} = \left(\frac{1}{2^{2\g} - 1}\right)^{1/2}\quad \mbox{and}
\quad C'_{\g} = \frac{2^{-\g}}{1 - 2^{-\g}}.$$ Step 2 of the proof
of \cite[Theorem~1.1]{Ibrahim09} asserts that:
\begin{equation}\label{fer1}
A_{1} \leq C \|g\|_{L^{\infty}},
\end{equation}
while \cite[Lemma 3.1]{Ibrahim09} gives:
\begin{equation}\label{fer2}
A_{2} \leq C \|f\|_{BMO}.
\end{equation}
In order to estimate $A_{3}$, we proceed in the following way:
$$
A_{3} \leq \sup_{j\geq 1} 2^{\g j} \|\vphi_{j}
* f\|_{L^{\infty}}\leq \sup_{j\geq 1} 2^{\g j} \|\vphi_{j}
* f\|_{L^{\infty}} + \|\vphi_{0}
* f\|_{L^{\infty}},\quad \vphi_{0} \mbox{ is given by
(\ref{new_seq})},
$$
hence (see Definition~\ref{def2})
$$
A_{3} \leq \|f\|_{B^{\g}_{\infty,\infty}}.
$$
Using the well known result (see for instance \cite{Farkas})
$$B^{\g}_{\infty,\infty} = C^{\g,\g/2},$$
we finally obtain
\begin{equation}\label{fer3}
A_{3} \leq \|f\|_{C^{\g,\g/2}}.
\end{equation}
Inequalities (\ref{fer0}), (\ref{fer1}), (\ref{fer2}) and
(\ref{fer3}) imply:
$$
\|f\|_{L^{\infty}}\leq C \left((2N+1)^{1/2} \|f\|_{BMO} + 2^{-\g N}
(\|f\|_{C^{\g,\g/2}}+ \|g\|_{L^{\infty}})\right).
$$
Optimizing the above inequality with respect to the variable $N$
(see Step 2 of the proof of \cite[Lemma~3.2]{Ibrahim09}), we
directly arrive into the result. $\hfill{\Box}$\\

We now present the proof of Theorem~\ref{theo2} that involve finer
estimates on the H\"{o}lder norm.\\

\noindent \textbf{Proof of Theorem~\ref{theo2}.} For the sake of
simplifying the ideas of the proof, we only consider $1$-spatial
dimensions $x=x_{1}$. The general $n$-dimensional case can be easily
deduced. Following the same notations of \cite{Ibrahim09}, we let
$\widetilde{\O}_{T} = (-1,2) \times (-T,2T)$, $\mathcal{Z}_{1}
\subseteq \mathcal{Z}_{2} \subseteq \widetilde{\O}_{T}$ such that
$$\mathcal{Z}_{1} = \{(x,t);\, -1/4<x<5/4 \,\mbox{ and } -T/4<t<5T/4\}$$
and
$$\mathcal{Z}_{2} = \{(x,t);\, -3/4<x<7/4 \,\mbox{ and } -3T/4<t<7T/4\}.$$
We also take the cut-off function $\Psi\in C^{\infty}_{0}(\R^{2})$,
$0\leq\Psi\leq 1$ satisfying:
\begin{equation}\label{5sar}
\Psi(x,t) = \left\{
\begin{aligned}
& 1 \quad \mbox{for}\quad (x,t)\in \mathcal{Z}_{1}\\
& 0 \quad \mbox{for}\quad (x,t)\in \R^{2}\setminus \mathcal{Z}_{2}.
\end{aligned}
\right.
\end{equation}
The main idea of the proof consists in extending the function $f$ to
a suitable function of the form $ \Psi\tilde{f}$ where $\tilde{f}$
is defined on $\widetilde{\O}_{T}$. We then apply inequality
(\ref{Ib:eq4}) (the scalar-valued version with $n=1$) to $
\Psi\tilde{f}$ and we estimate the different norms in order to get
the result. However, away from the complicated extension (Sobolev
extension) of the function $\tilde{f}$ that was done in
\cite{Ibrahim09}, we here consider a simpler symmetric extension.
Indeed, we first take the spatial symmetry of the function $f$:
\begin{equation}\label{spa_sym}
\tilde{f}(x,t) = \left\{
\begin{aligned}
& f(-x,t) \quad &\mbox{for}&\quad &-1<x<0&,\, &0\leq t\leq T&\\
& f(2-x,t) \quad &\mbox{for}&\quad &1<x<2&,\, &0\leq t\leq T&,
\end{aligned}
\right.
\end{equation}
and then the symmetry with respect to $t$:
\begin{equation}\label{tme_sym}
\tilde{f}(x,t) = \left\{
\begin{aligned}
& f(x,-t) \quad &\mbox{for}&\quad &-1<x<2&,\, &-T<t\leq 0& \\
& f(x,2T-t) \quad &\mbox{for}&\quad &-1<x<2&,\, &T\leq t < 2T&.
\end{aligned}
\right.
\end{equation}
We claim that $\Psi \tilde{f} \in C^{\g,\g/2}(\R^{2})$ with
\begin{equation}\label{h}
\|\Psi \tilde{f}\|_{C^{\g,\g/2}(\R^{2})} \leq
\|f\|_{C^{\g,\g/2}(\O_{T})}.
\end{equation}
In this case, we apply the scalar-valued version of inequality
(\ref{Ib:eq4}) (see Remark~\ref{rem1}) to the function $\Psi
\tilde{f}$ with $i=1$ and $g(x,t) =
\int_{0}^{x} \Psi(y,t) \tilde{f}(y,t) dy$. This, together with the
fact that $\Psi = 1$ on $\O_T$, lead to the following estimate:
\begin{equation}\label{estimate1}
\|f\|_{L^{\infty}(\O_{T})}\leq \|\Psi
\tilde{f}\|_{L^{\infty}(\R^{2})} \leq C \left(1 + \|\Psi
\tilde{f}\|_{BMO(\R^{2})}\left(\log^{+}(\|\Psi \tilde{f}\|_{
C^{\g,\g/2}(\R^{2})} +
\|g\|_{L^{\infty}(\R^{2})})\right)^{1/2}\right).
\end{equation}
It is worth noticing that choosing $i=1$ above is somehow
restrictive. In fact, we could also have used the inequality with
$i=2$ and $g(x,t) = \int_{0}^{t} \Psi(x,s) \tilde{f}(x,s) ds$.

\begin{figure}[!h]
\psfrag{Ott}{\scriptsize{$\widetilde{\O}_{T}$}}
\psfrag{0}{\scriptsize{$0$}}
\psfrag{T}{\scriptsize{$T$}}
\psfrag{Ot}{\scriptsize{$\O_{T}$}}
\psfrag{x}{\scriptsize{$x'$}}
\psfrag{xp}{\scriptsize{$x$}}
\psfrag{l}{\scriptsize{$\bar{x}$}}
%\psfrag{lp}{\scriptsize{$\bar{x}$}}
\psfrag{t}{\scriptsize{$t$}}
\psfrag{tb}{\scriptsize{$\bar{t}$}}
\begin{center}
\epsfig{file=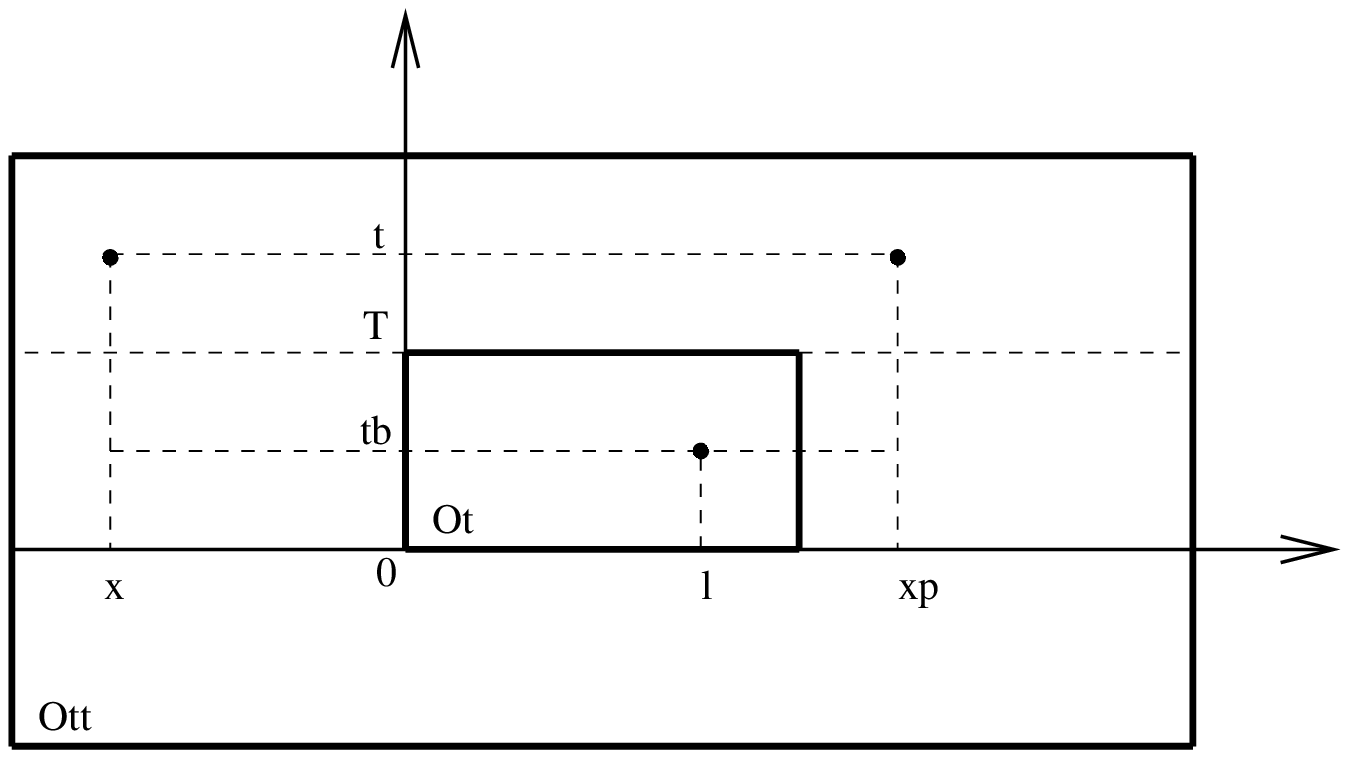, width=80mm}
\end{center}
\caption{Case $(x,t), (x',t)\in \widetilde{\O}_{T}$ with $\bar{x} = \bar{x}'$.}\label{fig2}
\end{figure}

In \cite{IM09} it was shown that $\|\Psi \tilde{f}\|_{BMO(\R^{2})}
\leq C (\|f\|_{BMO(\O_{T})} + \|f\|_{L^{1}(\O_{T})})$, while it is
clear that $\|g\|_{L^{\infty}(\R^{2})} \leq C
\|\tilde{f}\|_{L^{\infty}(\widetilde{\O}_T)}\leq C
\|f\|_{C^{\g,\g/2}(\O_{T})}$. These arguments, along with (\ref{h})
and (\ref{estimate1}), directly terminate the proof. The only point
left is to show the claim (\ref{h}). Recall the norm
$$\|\Psi \tilde{f}\|_{C^{\g,\g/2}(\R^{2})} = \|\Psi
\tilde{f}\|_{L^{\infty}(\R^{2})} +
\<\Psi\tilde{f}\>^{(\g)}_{x,\R^{2}} +
\<\Psi\tilde{f}\>^{(\g/2)}_{t,\R^{2}}.$$ It is evident that
$$\|\Psi \tilde{f}\|_{L^{\infty}(\R^{2})} \leq C
\|f\|_{L^{\infty}(\O_{T})},$$ hence we only need to estimate the two
terms $\<\Psi\tilde{f}\>^{(\g)}_{x,\R^{2}}$ and
$\<\Psi\tilde{f}\>^{(\g/2)}_{t,\R^{2}}$. We only deal with
$\<\Psi\tilde{f}\>^{(\g)}_{x,\R^{2}}$ since the second term can be
treated similarly. We examine the different positions of $(x,t),
(x',t)\in \R^{2}$. If $(x,t), (x',t)\in \R^{2}\setminus
\mathcal{Z}_{2}$, $x\neq x'$, then (since $\Psi=0$ over
$\R^{2}\setminus \mathcal{Z}_{2}$):
\begin{equation}\label{7alla1}
\frac{|(\Psi\tilde{f})(x,t) - (\Psi\tilde{f})(x',t)|}{|x-x'|^{\g}} =
0.
\end{equation}
If both $(x,t), (x',t)\in \widetilde{\O}_{T}$, $x\neq x'$, then the
special extension (\ref{spa_sym}) and (\ref{tme_sym}) of the
function $f$ guarantees the existence of
$$(\bar{x},\bar{t}), (\bar{x}',\bar{t})\in \O_T$$
such that:
\begin{equation}\label{D3T-1}
\tilde{f}(x,t) = f(\bar{x},\bar{t}),\quad \tilde{f}(x',t) =
f(\bar{x}',\bar{t}).
\end{equation}
Two cases can be considered. Either $\bar{x} = \bar{x}'$ (see Figure~\ref{fig2}), then we
forcedly have
$$\tilde{f}(x,t) = \tilde{f}(x',t),$$
\noindent and therefore
\begin{eqnarray}\label{Br}
\frac{|(\Psi\tilde{f})(x,t) - (\Psi\tilde{f})(x',t)|}{|x-x'|^{\g}}
&\leq&
\<\Psi\>^{(\g)}_{x,\widetilde{\O}_T}\|\tilde{f}\|_{L^{\infty}(\widetilde{\O}_{T})}\nonumber\\
&\leq& C \|f\|_{L^{\infty}(\O_{T})}\leq C
\|f\|_{C^{\g,\g/2}(\O_{T})},
\end{eqnarray}
or $\bar{x} \neq \bar{x}'$, then we forcedly have (see Figure~\ref{fig1})
\begin{equation}\label{D3T-0}
|x-x'|^{\g} \geq |\bar{x} - \bar{x}'|^{\g}.
\end{equation}
In this case, we compute:
\begin{figure}[!h]
\psfrag{Ott}{\scriptsize{$\widetilde{\O}_{T}$}}
\psfrag{0}{\scriptsize{$0$}}
\psfrag{T}{\scriptsize{$T$}}
\psfrag{Ot}{\scriptsize{$\O_{T}$}}
\psfrag{x}{\scriptsize{$x'$}}
\psfrag{xp}{\scriptsize{$x$}}
\psfrag{l}{\scriptsize{$\bar{x}'$}}
\psfrag{lp}{\scriptsize{$\bar{x}$}}
\psfrag{t}{\scriptsize{$t$}}
\psfrag{tb}{\scriptsize{$\bar{t}$}}
\begin{center}
\epsfig{file=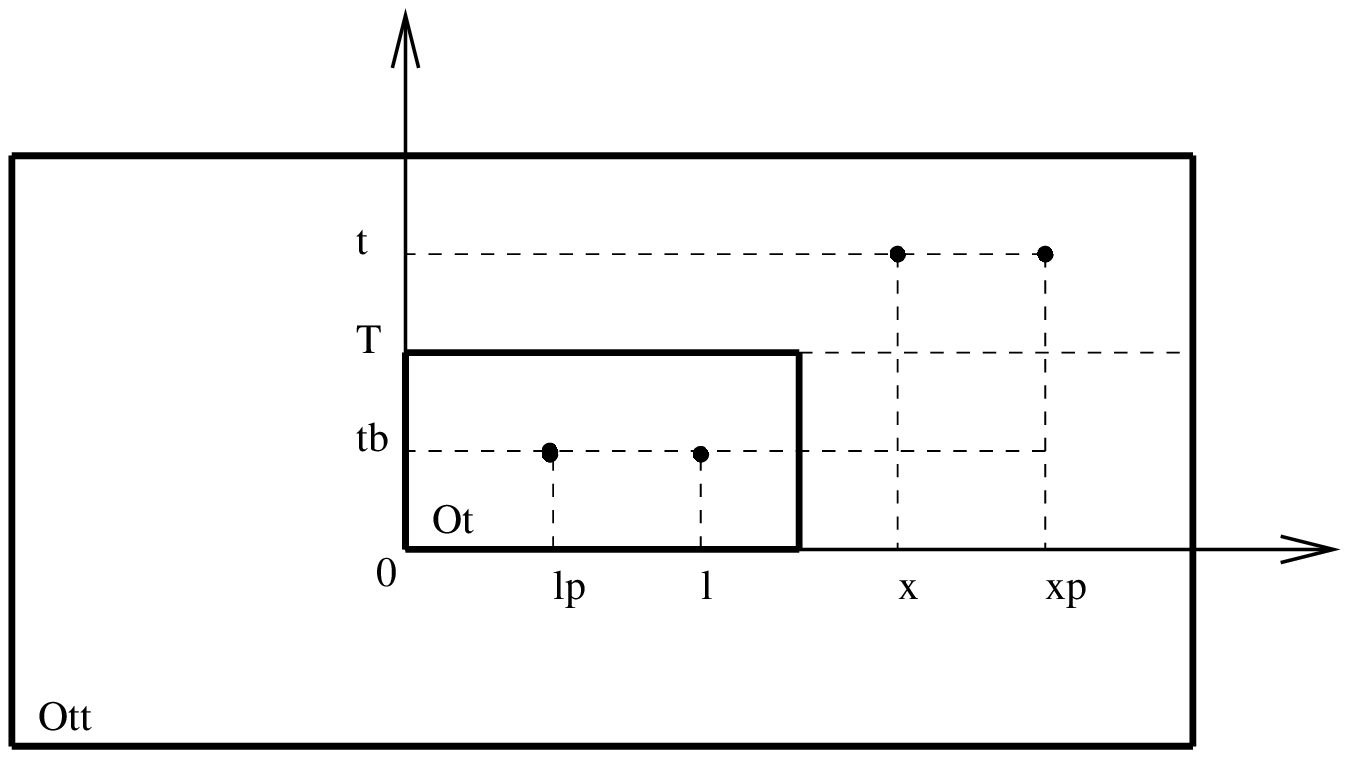, width=80mm} \hspace{0.2cm}
\epsfig{file=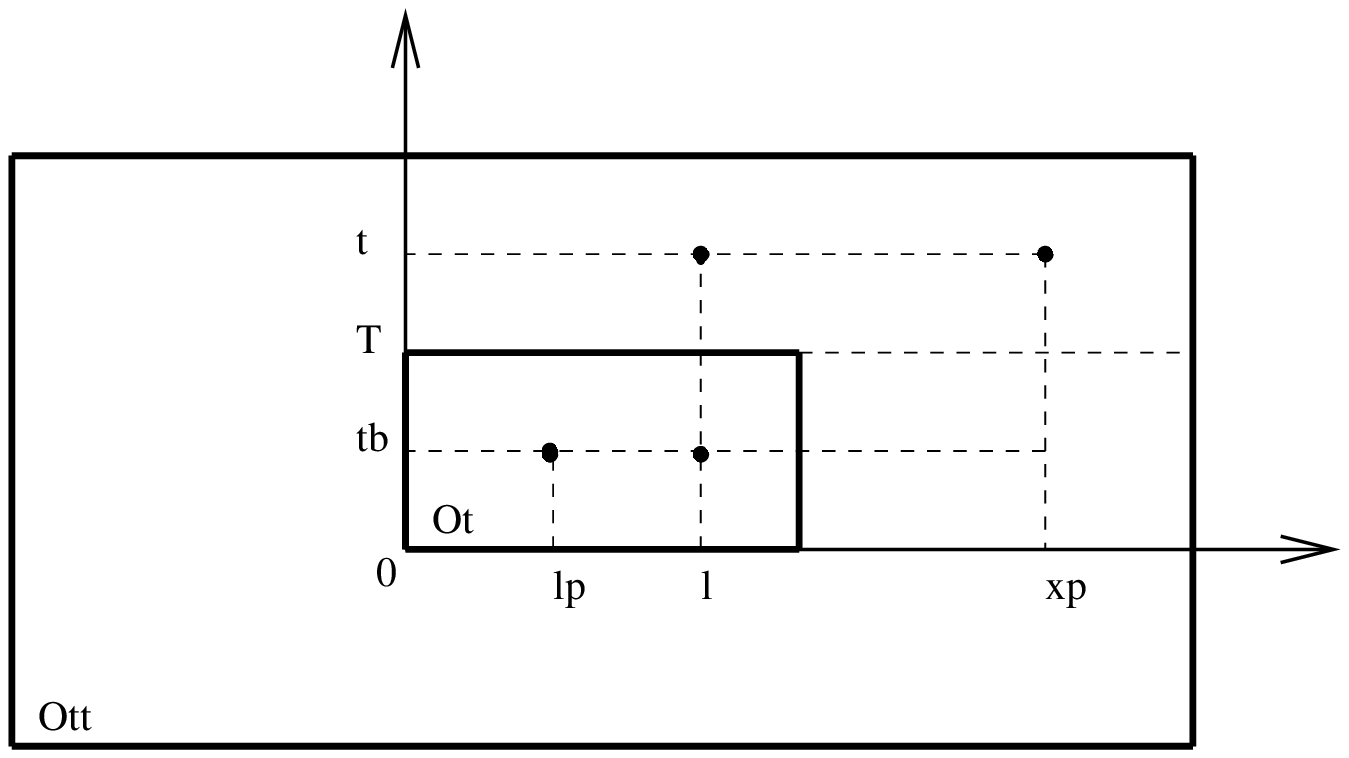, width=80mm}
\end{center}
\caption{Case $(x,t), (x',t)\in \widetilde{\O}_{T}$ with $\bar{x}\neq
  \bar{x}'$. On the right: $x' = \bar{x}'$. On the left: $x' \neq
  \bar{x}'$.}\label{fig1}
\end{figure}
\begin{eqnarray}\label{D3T1}
\frac{|(\Psi\tilde{f})(x,t) - (\Psi\tilde{f})(x',t)|}{|x-x'|^{\g}}
&\leq&  \frac{|\tilde{f}(x,t)||\Psi(x,t) - \Psi(x',t)|}{|x-x'|^{\g}}
+ \frac{|\Psi(x',t)||\tilde{f}(x,t) -
\tilde{f}(x',t)|}{|x-x'|^{\g}}\nonumber\\
&\leq&
\|\tilde{f}\|_{L^{\infty}(\widetilde{\O}_{T})}\<\Psi\>^{(\g)}_{x,\widetilde{\O}_T}
+\frac{|\tilde{f}(x,t) - \tilde{f}(x',t)|}{|x-x'|^{\g}}.
\end{eqnarray}
Using (\ref{D3T-1}) and (\ref{D3T-0}), we deduce that:
$$\frac{|\tilde{f}(x,t) - \tilde{f}(x',t)|}{|x-x'|^{\g}} =
\frac{|f(\bar{x},\bar{t}) - f(\bar{x}',\bar{t})|}{|x-x'|^{\g}}\leq
\frac{|f(\bar{x},\bar{t}) - f(\bar{x}',\bar{t})|}{|\bar{x} -
\bar{x}'|^{\g}}\leq \<f\>^{(\g)}_{x,\O_T},$$ therefore, by
(\ref{D3T1}), we obtain:
\begin{equation}\label{7alla2}
\frac{|(\Psi\tilde{f})(x,t) - (\Psi\tilde{f})(x',t)|}{|x-x'|^{\g}}
\leq
\|\tilde{f}\|_{L^{\infty}(\widetilde{\O}_{T})}\<\Psi\>^{(\g)}_{x,\widetilde{\O}_T}
+ \<f\>^{(\g)}_{x,\O_T}\leq C \|f\|_{C^{\g,\g/2}(\O_{T})}.
\end{equation}
The remaining case is when $(x,t)\in \mathcal{Z}_{2}$ and $(x',t)\in
\R^{2}\setminus \widetilde{\O}_T$ (see Figure~\ref{fig3}). In this case, we have
$(\Psi\tilde{f})(x',t) = 0$ and
\begin{equation}\label{fin_add}
|x-x'|^{\g} \geq \left(\frac{1}{4}\right)^{\g},
\end{equation}
hence
\begin{equation}\label{7alla3}
\frac{|(\Psi\tilde{f})(x,t) -
(\Psi\tilde{f})(x',t)|}{|x-x'|^{\g}}\leq 4^{\g}
\|\tilde{f}\|_{L^{\infty}(\mathcal{Z}_{2})}\leq C
\|f\|_{C^{\g,\g/2}(\O_{T})}.
\end{equation}
\begin{figure}[!h]
\psfrag{Ott}{\scriptsize{$\widetilde{\O}_{T}$}}
\psfrag{0}{\scriptsize{$0$}}
\psfrag{T}{\scriptsize{$T$}}
\psfrag{Ot}{\scriptsize{$\O_{T}$}}
\psfrag{Z2}{\scriptsize{$\mathcal{Z}_{2}$}}
\psfrag{x}{\scriptsize{$x$}}
\psfrag{xp}{\scriptsize{$x'$}}
\psfrag{t}{\scriptsize{$t$}}
\psfrag{1}{\tiny{\hspace{-0.2cm}{$1/4$}}}
\psfrag{To4}{\tiny{$T/4$}}
\begin{center}
\epsfig{file=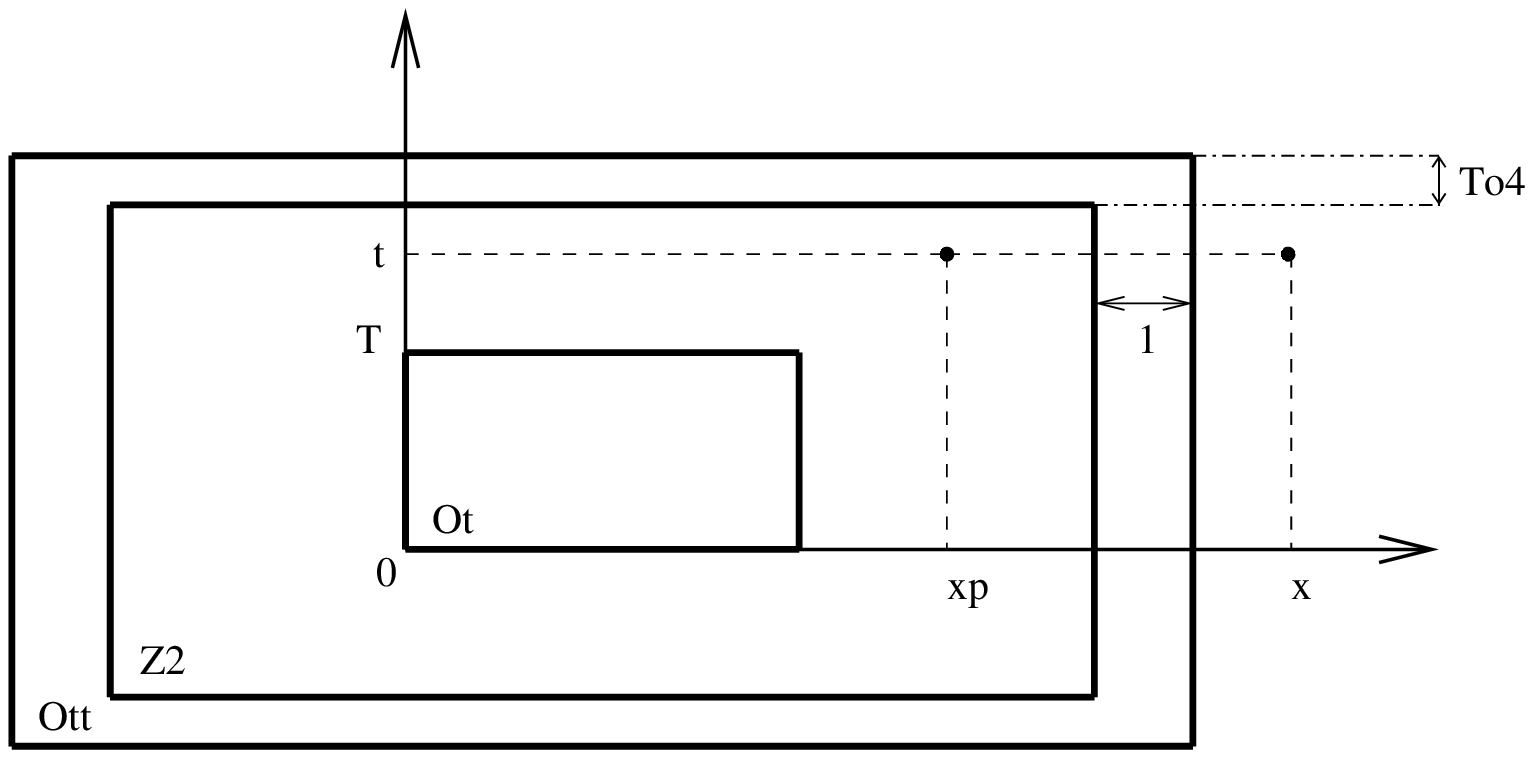, width=85mm}
\end{center}
\caption{case $(x,t)\in \mathcal{Z}_{2}$ and $(x',t)\in
\R^{2}\setminus \widetilde{\O}_T$.}\label{fig3}
\end{figure}
From (\ref{7alla1}), (\ref{Br}), (\ref{7alla2}) and (\ref{7alla3}),
we finally deduce that
$$
\<\Psi\tilde{f}\>^{(\g)}_{x,\R^{2}} \leq C
\|f\|_{C^{\g,\g/2}(\O_T)}.
$$
Arguing in exactly the same way as above, we also find that:
$$
\<\Psi\tilde{f}\>^{(\g/2)}_{t,\R^{2}} \leq C
\|f\|_{C^{\g,\g/2}(\O_T)},
$$
with a possibly different constant $C$ that depend on $T$. Indeed,
the term $T$ enters in estimating
$\<\Psi\tilde{f}\>^{(\g/2)}_{t,\R^{2}}$ since (\ref{fin_add}) is now
replaced (see again Figure~\ref{fig3}) by
$$|t-t'|^{\g} \geq \left(\frac{T}{4}\right)^{\g}.$$
This shows the claim. $\hfill{\Box}$
\begin{rem}
In the case of multi-spatial coordinates $x_{i}$, $i = 1,\ldots,n$,
we simultaneously apply the extension (\ref{spa_sym}) to each spatial
coordinate while fixing all other coordinates including $t$. Finally,
fixing the spatial variables, we make the extension with respect to
$t$ as in (\ref{tme_sym}).
\end{rem}

%\bibliographystyle{siam}
%\bibliography{biblio}

\begin{thebibliography}{20}
\setlength{\itemsep}{5pt}

\bibitem{BG80} H. Br\'{e}zis and T. Gallou\"{e}t, Nonlinear
Schr\"{o}dinger evolution equations, Nonlinear Anal., 4 (1980), pp.
677-681.


\bibitem{BW80} H. Br\'{e}zis and S. Wainger, A note on limiting
cases of Sobolev embeddings and convolution inequalities, Comm.
Partial Differential Equations, 5 (1980), pp. 773-789.


\bibitem{Engler} H. Engler, An alternative proof of the
Brezis-Wainger inequality, Comm. Partial Differential Equations, 14
(1989), pp. 541-544.


\bibitem{Farkas} W. Farkas, J. Johnsen and W. Sickel, Traces of anisotropic
Besov-Lizorkin-Triebel spaces---a complete treatment of the
borderline cases, Math. Bohem., 125 (2000), pp. 1-37.


\bibitem{Ibrahim09} H. Ibrahim, A critical parabolic Sobolev embedding via
Littlewood-Paley decomposition, preprint, arxiv:0908.1866v1.






\bibitem{IJM08} H. Ibrahim, M. Jazar and R. Monneau, Global
existence of solutions to a singular parabolic/Hamilton-Jacobi
coupled system with Dirichlet conditions, C. R. Math. Acad. Sci.
Paris, 346 (2008), pp. 945-950.





\bibitem{IM09} H. Ibrahim and R. Monneau, On a parabolic logarithmic
sobolev inequality, J. Funct. Anal., 257 (2009), pp. 903-930.






\bibitem{KOT02} H. Kozono, T. Ogawa and Y. Taniuchi, The critical
Sobolev inequalities in Besov spaces and regularity criterion to
some semi-linear evolution equations, Math. Z., 242 (2002), pp.
251-278.

\bibitem{KOT03} H. Kozono, T. Ogawa and Y. Taniuchi, Navier-Stokes equations
in the Besov space near $L\sp\infty$ and BMO, Kyushu J. Math., 57
(2003), pp. 303-324.




\bibitem{KT00} H. Kozono and Y. Taniuchi, Limiting case of the
Sobolev inequality in $BMO$, with application to the Euler
equations, Comm. Math. Phys., 214 (2000), pp. 191-200.


\bibitem{LSU} O. A. Lady{\v{z}}enskaja, V. A. Solonnikov and N. N.
Ural$'$ceva, Linear and quasilinear equations of parabolic type,
Translated from the Russian by S. Smith. Translations of
Mathematical Monographs, Vol. 23, American Mathematical Society,
Providence, R.I., 1967.

\bibitem{Og03} T. Ogawa, Sharp Sobolev inequality of logarithmic
type and the limiting regularity condition to the harmonic heat
flow, SIAM J. Math. Anal., 34 (2003), pp. 1318-1330 (electronic).



\bibitem{Tri} H. Triebel, Theory of function spaces, vol. 78 of
Monographs in Mathematics, Birkh\"{a}user Verlag, Basel, 1983.


\bibitem{Tri1} H. Triebel, Theory of function spaces. $III$, vol.
100 of Monographs in Mathematics, Birkh\"{a}user Verlag, Basel,
2006.





\end{thebibliography}
\end{document}